\documentclass[12pt]{amsart}

\usepackage{a4wide}
\usepackage{amssymb}
\usepackage{amsfonts}
\usepackage{amsmath}
\usepackage[margin=3cm]{geometry}

\date{}

\newtheorem{proposition}{Proposition}[section]
\newtheorem{theorem}[proposition]{Theorem}
\newtheorem{lemma}[proposition]{Lemma}

\newtheorem{corollary}[proposition]{Corollary}

\def\GK{{\rm  GK}\,}

\def\der{\partial }

\def\nFM0{{\nu }_{F,M_0}}
\def\nFN0{{\nu }_{F,N_0}}
\def\nGN0{{\nu }_{G,N_0}}

\def\Nat{\mathbb{N}}

\def\N0{ {\bf N}_0 }

\def\ra{\rightarrow}

\def\Xpm{X^{\pm }}

\def\s{\sigma}
\def\Z{\mathbb{Z}}

\def\l1{{\lambda}_1}

\def\a{\alpha}
\def\a0{ {\alpha }_0}
\def\a1{ {\alpha }_1}

\def\l{\lambda}

\def\nFGM0{{\nu }_{F,G,M_0}}

\def\nFN0{{\nu}_{F,N_0}}

\def\sm{{\sigma}^m}

\def\sm1{{\sigma}^{-1}}

\def\smtp1{{\sigma}^{-t+1}}

\def\S1{S^{-1}}

\def\Xpm1{X^{\pm 1}_1}

\def\sPM1{{\sigma }^{\pm 1}}
\def\sMP1{{\sigma }^{\mp 1 }}

\def\Ytm1{Y^{t-1}}
\def\Yim1{Y^{i-1}}

\def\Aut{{\rm Aut}}

\def\CI{\mathbb{I}}

\def\marginpar{}

\begin{document}

\author{V. V.  Bavula}
\address{Department of Pure Mathematics \\ University of Sheffield \\ Hicks Building \\ Sheffield S3  7RH, UK}
\email{v.bavula@sheffield.ac.uk}

\author{V. Levandovskyy }
\address{Lehrstuhl D f\"ur Mathematik \\ RWTH Aachen University \\ 52062 Aachen, Germany} \email{Viktor.Levandovskyy@math.rwth-aachen.de}


\title{A remark on the Dixmier Conjecture}

\maketitle
\begin{abstract}
The Dixmier Conjecture says that every endomorphism of the
(first) Weyl algebra $A_1$ (over a field of characteristic zero)
is an automorphism, i.e.,
if $PQ-QP=1$ for some $P, Q \in A_1$ then $A_1 = K \langle P, Q \rangle$.
The Weyl algebra $A_1$ is a $\Z$-graded algebra. We prove that the
Dixmier Conjecture holds if the elements $P$ and $Q$ are sums of no more than
two homogeneous elements of $A$ (there is no restriction on the total degrees of $P$ and $Q$).\\

{\em Key Words: the Weyl algebra, the Dixmier Conjecture, automorphism, endomorphism, a $\Z$-graded algebra.}

 {\em Mathematics subject classification 2010: 16S50, 16W20, 16S32, 16W50.}
\end{abstract}

\section{Introduction}

In the paper, $K$ is a field of characteristic zero and $K^*:=K \setminus\{0\}$. The algebra
$A_1 := K \langle X, Y \mid [Y, X]=1 \rangle$ is called the {\em first Weyl algebra} where
$[Y, X]= YX - XY$. The $n$'th tensor power of $A_1$, $A_n := A_1^{\otimes n} =
\underbrace{A_1 \otimes \cdots \otimes A_1}_{n \text{ times}}$, is called the {\em $n$'th Weyl algebra}.
The algebra $A_n$ is a simple Noetherian domain of Gel'fand-Kirillov dimension $\GK(A_n)=2n$,
it is canonically isomorphic to the algebra of polynomial differential operators
$K \langle X_1, \ldots, X_n, \partial_1, \ldots, \partial_n \rangle$ (where $\partial_i = \tfrac{\partial}{\partial x_i}$)
via $X_i \mapsto X_i$, $Y_i \mapsto \partial_i$ for $i=1,\ldots,n$.

In his seminal paper \cite{Dix}, Dixmier (1968) found explicit generators for the group $G = \Aut_K(A_1)$ of
$K$-automorphisms of the Weyl algebra $A_1$. Namely, the group $G$ is generated by the
obvious automorphisms:
\[
(X, Y) \mapsto (X, Y + \l X^n), \quad (X, Y) \mapsto (X + \l Y^n, Y), \quad  (X, Y) \mapsto (\mu X, \mu^{-1} Y)
\]
where $\l \in K$, $\mu \in K^*$ and $n \in \Nat_{+}:=\{1, 2, \ldots\}$.

In \cite{Dix}, Dixmier posed six problems: The first problem of Dixmier  (in the list) asks if {\em every endomorphism of the Weyl algebra $A_1$ is an automorphism}, i.e., given elements $P, Q$ of $A$ such that $[P, Q]=1$, do
they generate the algebra $A_1$? A similar problem but for the $n$'th Weyl algebra is called the {\em Dixmier
Conjecture}. Problems 3 and 6 have been solved by Joseph \cite{josclA1} (1975), Problem 5 and Problem 4 (in the case
of homogeneous elements) have been solved by Bavula \cite{Bav-DixPr5} (2005).

The Dixmier Conjecture implies the {\em Jacobian Conjecture} (see \cite{BCW}) and the inverse implication
is also true (see \cite{Tsuchimoto-2005} and \cite{Belov-Konts-2007}); a short proof is given in \cite{Bav-DCJC-2005}; see also \cite{Adjam-vdEs-2007}).

In \cite{Bav-RenQn}, it is shown that for each $K$-endomorphism $\phi: A_n \to A_n$ its image is very large, i.e.,
the left $A_{2n}$-module ${}^{\phi}{A_n}^{\phi}$ is a holonomic $A_{2n}$-module (where
for all $a, b \in A_n$ and $c \in {}^{\phi}{A_n}^{\phi}$, $a \cdot c \cdot b := \phi(a) c \phi(b)$). In particular, it has finite length
 with simple holonomic factors over $A_{2n}$ (see \cite{Bav-RenQn} for details). To prove that the Dixmier Conjecture holds for the Weyl algebra $A_n$ it remains to show that the length is 1. Note, that the Gel'fand-Kirillov dimension of a simple $A_{2n}$-module can be $2n, 2n+1, \ldots, 4n-1$,  and the last case is the generic case.

In \cite{Bav-T1DC}, it is shown that every algebra endomorphism of the algebra $\CI_1 = K\langle x, \partial, \int \rangle$
of polynomial integro-differential operators is an automorphism and it is conjectured that the same result holds for $\CI_n := \CI_1^{\otimes n}=K\langle x_1, \ldots , x_n, \der_1, \ldots , \der_n , \int_1, \ldots , \int_n\rangle$.

The Weyl algebra $A = \oplus_{i\in\Z} A_{1,i}$ is a $\Z$-graded algebra ($A_{1,i} A_{1,j} \subseteq A_{1,i+j}$ for all
$i, j\in \Z$) where $A_{1,0}=K[H]$, $H=YX$ and, for $i\geq 1$, $A_{1,i} = K[H]X^i$ and $A_{1,-i} = K[H]Y^i$.
For a nonzero element $a$ of $A_1$, the number of {\em nonzero homogeneous} components is called the {\em mass} of $a$,
denoted by $m(a)$. For example, $m(\alpha X^i)=1$ for all $\alpha \in K[H]\setminus\{0\}$ and $i\geq 1$. The aim of this paper is to prove the following theorem.

\begin{theorem}\label{20Sep16} 
Let $P, Q$ be elements of the first Weyl algebra $A_1$ with $m(P)\leq 2$ and $m(Q)\leq 2$.
If $[P,Q]=1$ then $P = \tau(Y)$ and $Q=\tau(X)$ for some automorphism $\tau\in\Aut_K (A_1)$.
\end{theorem}

\section{Proof of Theorem \ref{20Sep16}}


{\bf The Weyl algebra is a generalized Weyl algebra}.
Let $D$ be a ring with an automorphism $\s $ and a central element $a$.
 The {\bf generalized Weyl algebra} $A=D(\sigma, a)$  of degree 1,
is the ring generated by $D$ and two indeterminates $X$ an $Y$  subject to
the relations \cite{Bav-FA-1991}:
$$
X\alpha=\sigma(\alpha)X \ {\rm and}\ Y\alpha=\sigma^{-1}(\alpha)Y,\; {\rm for \; all }\;
\alpha \in D, \ YX=a \ {\rm and}\ XY=\sigma(a).
 $$
The algebra $A={\oplus}_{n\in \mathbb{Z}}\, A_n$ is a
$\mathbb{Z}$-graded algebra where $A_n=Dv_n$,
 $v_n=X^n\,\, (n>0), \,\,v_n=Y^{-n}\,\, (n<0), \,\,v_0=1.$
 It follows from the defining relations that
$$v_nv_m=(n,m)v_{n+m}=v_{n+m}<n,m> $$
for some elements $(n,m)=\s^{-n-m}(<n,m>)\in D$. If $n>0$ and $m>0$ then
\begin{eqnarray*}
n\geq m & : & (n,-m)=\sigma^n(a)\cdots \sigma^{n-m+1}(a), \,\,  (-n,m)=\sigma^{-n+1}(a)\cdots \sigma^{-n+m}(a), \\
n\leq m & : & (n,-m)=\sigma^{n}(a)\cdots \sigma(a),\,\,\, (-n,m)=\sigma^{-n+1}(a)\cdots a,
\end{eqnarray*}
in other cases $(n,m)=1$.

Let $K[H]$ be a polynomial ring in a variable $H$ over the field $K$, $\s :H\ra H-1$ be the
$K$-automorphism of the algebra $K[H]$ and $a=H$. The first Weyl algebra
$A_1=K\langle X, Y \mid YX-XY=1\rangle$ is isomorphic to the generalized Weyl algebra
$$A_1\simeq K[H](\s , H),\; X \mapsto X,\; Y \mapsto Y,\; YX  \mapsto H.$$
 We identify both these algebras via this isomorphism, that is
$A_1=K[H](\s , H)$ and $H=YX$.

If $n>0$ and $m>0$ then
\begin{eqnarray*}
n\geq m & : & (n,-m)= (H-n)\cdots (H-n+m-1),\; (-n,m)=(H+n-1)\cdots (H+n-m), \\
n\leq m & : & (n,-m)=(H-n)\cdots (H-1), \; (-n,m)=(H+n-1)\cdots H,
\end{eqnarray*}
in other cases $(n,m)=1$.

The localization $B=S^{-1}A_1$ of the Weyl algebra $A_1$ at the Ore subset
 $S=K[H]\backslash \{ 0\}$ of $A_1$ is
the {\em skew Laurent polynomial ring} $B=K(H)[X, X^{-1}; \s ]$ with
coefficients from the field $K(H)=S^{-1}K[H]$ of rational functions
where $\s \in \Aut_K \, K(H)$ and $\s (H)=H-1$. The map $A_1\ra B$, $a  \mapsto a/1$
is an algebra monomorphism. We identify the algebra $A_1$ with its image in the algebra $B$
via $A_1 \ra  B, \;\; X \mapsto  X, \;\; Y \mapsto HX^{-1}$.
 The algebra $B=\oplus_{i\in \mathbb{Z}}\, B_i$ is a $\mathbb{Z}$-graded algebra where
$B_i=K(H)X^i$. The algebra $A_1$ is a $\mathbb{Z}$-graded subalgebra of $B$.

A polynomial $f(H)=\l_nH^n+\l_{n-1}H^{n-1}+\cdots +\l_0\in K[H]$
of degree $n$ is called a {\em  monic}  polynomial  if the {\em leading coefficient}
  $\l_n$ of $f(H)$ is $1$. A rational function $h\in
K(H)$ is called a {\em monic}  rational function  if $h=f/g$ for
some monic polynomials $f,g$. A homogeneous element $u=\alpha x^n$
of $B$ is called {\em monic} if $\alpha $ is a monic rational function.
We can extend the concept of degree of polynomial to the field of rational functions by the rule
 $\deg \, h = \deg\, f - \deg\, g$ where $h=f/g\in K[H]$. If $h_1, h_2 \in K(H)$
then $\deg\, h_1h_2=\deg\,h_1 +\deg\, h_2$ and
 $\deg (h_1+h_2)\leq \max \{ \deg\,h_1 , \deg\, h_2\}$.
  We denote by ${\rm sign} (n)$ and by $|n|$ the {\em sign} and the
{\em absolute value} of $n\in \mathbb{Z}$, respectively.

Let $A$ be an algebra and $a\in A$. The subalgebra of $A$, $C_A(a)=\{ b \in A \mid ab=ba \}$, is called
the {\em centralizer} of the element $a$ in $A$.

\begin{proposition}[\cite{Bav-DixPr5}, Proposition 2.1]
\label{cheab} \marginpar{cheab}
{\em (Centralizer of a Homogeneous Element of the Algebra $B$)}
\begin{enumerate}
\item Let $u=\alpha X^n$ be a monic element of $B_n$  with $n\neq 0$.  Then the centralizer $C_B(u)=K[v,v^{-1}]$ is a Laurent polynomial
ring
for a unique element $v=\beta X^{{\rm sign} (n) s}$
where $s$ is the least positive divisor of $n$ for which there
exists an element $\beta =\beta_s \in K(H)$, necessarily monic and
uniquely defined, such that
\begin{equation}
\label{bap}
\beta \,\s^s (\beta )\, \s^{2s}(\beta ) \cdots \s^{(n/s -1)s}(\beta )=\alpha, \;\; {\rm if}\;\; n>0,
\end{equation}
\begin{equation}
\label{bam}
\beta \, \s^{-s} (\beta ) \,\s^{-2s}(\beta ) \cdots \s^{-(|n|/s -1)s}(\beta )=\alpha, \;\; {\rm if}\;\; n<0.
\end{equation}

\item Let $u\in K(H)\backslash K$. Then $C_B(u)=K(H)$.
\end{enumerate}
\end{proposition}

\noindent
Let $A_{1,+} := K[H][X; \sigma]$ and $A_{1,-} := K[H][Y; \sigma^{-1}]$. The algebras
$A_{1,+}$ and $A_{1,-}$ are (skew polynomial) subalgebras of $A_1$.
\begin{lemma}[\cite{Bav-DixPr5}]
\label{a20Sep16} \marginpar{a20Sep16}
If $u \in A_{1,\pm}\setminus\{0\}$ then $C_A(u) \subseteq A_{1,\pm}$.
\end{lemma}

\noindent
The $K$-automorphism of the Weyl algebra $A_1$,
\marginpar{xiaut}
\begin{equation}
\label{xiaut}
\xi: A_1 \to A_1, \; X \mapsto Y, \; Y \mapsto -X,
\end{equation}
reverses the $\Z$-grading of the Weyl algebra $A_1$, that is
\marginpar{xiaut1}
\begin{equation}
\label{xiaut1}
\xi(A_{1,i}) = A_{1,-i} \text{ for all } z\in\Z.
\end{equation}

\noindent
By the {\em degree} of an element of $A_1$ we mean its {\em total degree} with respect to
the canonical generators $X$ and $Y$ of $A_1$.
Let $A_{1, \leq i} := \{ p\in A \mid \deg(p)\leq i\}$ for $i\in\Nat$. Then
$\{A_{1, \leq i}\}_{i\in\Nat}$ is the standard filtration of the algebra $A_1$ associated with the generators $X$ and $Y$.
\noindent
For all $i\in\Z\setminus\{0\}$ and $f\in K[H]\setminus K$,
\marginpar{degsf}
\begin{equation}\label{degsf}
\deg \s^i(f) = \deg f \text{ and } \deg (1 - \s^i)(f)  = \deg f - 1.
\end{equation}

\noindent
{\bf Proof of Theorem \ref{20Sep16}}:
{\em (i)} {\em If $P, Q\in A_{1,\leq 1}$ then $P = \tau(Y)$ and $Q=\tau(X)$ for some $\tau\in\Aut_K (A_1)$}:
Clearly, $P = aY + bX + \l$ and $Q = cY +dX + \mu$ for some $a, b, c, d, \l, \mu \in K$. Then
$1 = [P,Q]=ad-bc$. So, the automorphism $\tau$ can be chosen of the form
$\tau(Y) = aY + bX + \l$ and $\tau(X) = cY +dX + \mu$.

So, till the end of the proof we assume that at least one of the polynomials $P$ or $Q$ does not belong to the
space $A_{1, \leq 1}$. In view of the relation $1 = [P,Q] = [-Q,P]$, we can assume that $P \notin A_{1, \leq 1}$.
In view of Equation (\ref{xiaut1}), we can assume that the highest homogeneous part of $P$, say $P_p \in A_{1,p}$,
satisfies the condition that $p \geq 2$. Since $m(P)\leq 2$, either $P = P_p$ (if $m(P)=1$) or otherwise
$P = P_r + P_p$ for some nonzero $P_r \in A_{1,r}$ where $r<p$.

{\em (ii)} $(m(P), m(Q)) \neq (1,1)$:
Suppose that $m(P)= m(Q)=1$, we seek a contradiction. Then $P =\alpha X^p$ and $Q=\beta Y^p$ for some nonzero polynomials
$\alpha, \beta\in K[H]$.  Then
\[
1 = [P,Q] = \alpha \s^p(\beta) (p,-p) - \beta \s^{-p}(\alpha) (-p,p)
\]
\[
=
\alpha \s^p(\beta) (p,-p) - \beta \s^{-p}(\alpha) \s^{-p}((p,-p)) =
(1-\s^{-p})(\alpha \s^p(\beta) (p,-p)).
\]
Since $p \geq 2$ (or $P\notin A_{1, \leq 1}$),
$0 = \deg 1 = \deg\, (1-\s^{-p})(\alpha \s^p(\beta) (p,-p)) =
\deg \alpha + \deg \beta +\deg\, (p,-p) -1$
(by Equation (\ref{degsf})) $\geq 0+0+p-1 \geq 2-1=1$, a contradiction.

{\em (iii)} $(m(P), m(Q)) \neq (1,2)$: Suppose that $m(P)=1$ and $m(Q)=2$. Then
$P = \alpha X^p$ for some $p\geq 2$ and $Q=Q_s + Q_q$ where $Q_s \in A_{1,s}$, $Q_q \in A_{1,q}$  and $s<q$.
By Lemma \ref{a20Sep16}, the equality $[P,Q]=1$ implies that $[P, Q_s]=1$ and $[P, Q_q]=0$. By the case (ii), this is not possible.

{\em (iv)} Suppose that $m(P)=2$ and $m(Q)=1$. Then $P = P_r + P_p$ and $Q=Q_q$.
By Lemma \ref{a20Sep16} the equality $[P,Q]=1$ implies that $[P_p, Q_q]=0$ and $[P_r, Q_q]=1$. Then, $q\geq 0$, by Lemma \ref{a20Sep16}.
The case $q=0$ is not possible since then both $P_r, Q_q\in K[H]$ and this would contradict  the equality $[P_r, Q_q]=1$.
Therefore, $q>0$.
Then $P_r = \beta Y^q$ and $Q_q = \alpha X^q$ for some nonzero elements
$\beta, \alpha \in K[H]$.  Then
\[
-1 = [Q_q, P_r] = (1-\s^{-q}) (\alpha \s^p(\beta) (q,-q) )
\]
implies that
$0 = \deg (-1) = \deg\, (1-\s^{-q}) (\alpha \s^p(\beta) (q,-q) ) = \deg \alpha + \deg \beta +q -1$,  by Equation (\ref{degsf}).
Hence, $q=1$, $\alpha, \beta \in K^*$ and $\beta = -\alpha^{-1}$. Then
$P, Q \in A_{1,\leq 1}$, and, by the statement (i), the pair $(P,Q)$ is obtained from the pair $(Y,X)$ by
applying an automorphism of $A_1$.

{\em (v)} $(m(P), m(Q)) \neq (2,2)$:
Since $m(P)=m(Q)=2$, we can write $P = P_r + P_p$ and $Q = Q_s + Q_q$ as sums
of homogeneous elements where $r < p$, $P_r \in A_{1,r}$, $P_p \in A_{1,p}$ and $s<q$,
$Q_s \in A_{1,s}$, $Q_q \in A_{1,q}$.
The equality $[P,Q]=1$ implies that $[P_r, Q_s]=0$ and $[P_p, Q_q]=0$ (see Lemma \ref{a20Sep16}). By Lemma \ref{a20Sep16}, the
elements $r$ and $s$ have the same sign (i.e., either $r<0, s<0$ or $r=s=0$ or $r>0, s>0$) and also
the elements $p$ and $q$ have the same sign.
Since $p\geq 2$, we must have $q>0$.

Suppose that $r\geq 0$, we seek a contradiction. Then $s\geq 0$ and
so the elements $P$ and $Q$ are elements of the subring $A_{1,+}= \oplus_{i \geq 0} K[H] X^i$.
Now,
$$
K[H] \ni 1 = [P, Q] \in [A_{1,+}, A_{1,+}] \subseteq \oplus_{i \geq 1} K[H] X^i,
$$ a contradiction.
Therefore, $r<0$ and $s<0$.

The equality $1=[P,Q]=[P_r, Q_q] + [P_p, Q_s]$ and Lemma \ref{a20Sep16} imply that $r+q=0$ and $p+s=0$,
that is $r=-q$ and $s=-p$. So,
$$
P = P_{-q} + P_p \text{ and } Q = Q_{-p} + Q_q.
$$
The elements $P_p$ and $P_{-q}$ are homogeneous elements of the Weyl algebra $A_1$.
The Weyl algebra $A_1$ is a homogeneous subalgebra of the algebra
$K(H)[X, X^{-1}; \s] = K(H)[Y, Y^{-1}; \s^{-1}]$ where $K(H)$ is the field of rational functions
in the variable $H$ and the automorphism $\s$ of $K(H)$ is given by the rule $\s(H) = H-1$.
By \cite[Proposition 2.1(1)]{Bav-DixPr5}, the centralizer $C_B(P_p)$ of the element $P_p$ in $B$ is a Laurent
polynomial algebra $K[\alpha X^n, (\alpha X^n)^{-1}]$ for some nonzero element
$\alpha\in K(H)$ and $n \geq 1$. In general, $\alpha \notin K[H]$.
Similarly, $C_B(P_{-q}) = K[\beta Y^m, (\beta Y^m)^{-1}]$ for some nonzero element
$\beta \in K(H)$ and $m \geq 1$.

Since $[P_p, Q_q]=0$, $Q_q \in C_B(P_p)$ and
\[
P_p = \lambda(P_p) (\alpha X^n)^i =
\lambda(P_p) \alpha  \s^n(\alpha) \cdots \s^{n(i-1)}(\alpha) X^{ni} = \alpha_{n,i} X^p,
\]
\[Q_q = \lambda(Q_q) (\alpha X^n)^j = \lambda(Q_q) \alpha \s^n(\alpha) \cdots \s^{n(j-1)}(\alpha) X^{nj} = \alpha'_{n,j} X^q,
\]
for some nonzero scalars $\lambda(P_p),  \lambda(Q_q) \in K^*$ and some $i \geq 1$ and $j\geq 1$
where
$$
\alpha_{n,i} = \lambda(P_p) \alpha  \s^n(\alpha) \cdots \s^{n(i-1)}(\alpha) \in K[H], \; p=ni,
$$
$$
\alpha'_{n,j}  = \lambda(Q_q) \alpha  \s^n(\alpha) \cdots \s^{n(j-1)}(\alpha) \in K[H], \; q=nj.
$$

Since $[P_{-p}, Q_{-p}]=0$, $Q_{-p} \in C_B(P_{-q})$ and
\[
P_{-q} = \lambda(P_{-q}) (\beta Y^m)^s =
\lambda(P_{-q}) \beta \s^{-m}(\beta)\cdots \s^{-m(s-1)}(\beta) Y^{ms} = \beta_{m,s} Y^p,
\]
\[Q_{-p} = \lambda(Q_{-p}) (\beta Y^m)^t =
\lambda(Q_{-p}) \beta\s^{-m}(\beta)\cdots \s^{-m(t-1)}(\beta) Y^{mt} = \beta'_{m,t} Y^q,
\]
for some nonzero scalars $\lambda(P_{-q}),  \lambda(Q_{-p}) \in K^*$ and some $s \geq 1$ and $t\geq 1$
where
$$
\beta_{m,s} = \lambda(P_{-q}) \beta \s^{-m}(\beta) \cdots \s^{-m(s-1)}(\beta) \in K[H],\;  p=ms,
$$
$$
\beta'_{m,t}  = \lambda(Q_{-p})\beta \s^{-m}(\beta) \cdots \s^{-m(t-1)}(\beta) \in K[H], \; q=mt.
$$

Now,
\[
1 = [P,Q] = [P_p, Q_{-p}] + [P_{-q}, Q_q]
 = [\alpha_{n,i} X^p, \beta_{m,t}' Y^p ] + [\beta_{m,s} Y^q , \alpha_{n,j}' X^q]
\]
\[
= \alpha_{n,i} \s^p (\beta'_{m,t}) (p,-p) - \beta'_{m,t} \s^{-p}(\alpha_{n,i}) (-p,p)
\]
\[+
\beta_{m,s} \s^{-q}(\alpha'_{n,j}) (-q,q) - \alpha'_{n,j} \s^{q}(\beta_{m,s}) (q,-q).
\]

Using the equalities $(-p,p)=\s^{-p}((p,-p))$ and $(-q,q)=\s^{-q}((q,-q))$, the last
equality above can be rewritten as follows
 \marginpar{1=ab}
\begin{equation}
\label{1=ab}
1 = (1-\s^{-p})(a) + (1-\s^{-q})(b)
\end{equation}
where $a = \alpha_{n,i} \s^{p}(\beta'_{m,t}) (p,-p) \in K[H]$ and $b = \alpha'_{n,j} \s^{q}(\beta_{m,s}) (q,-q) \in K[H]$.

Recall that $P = P_{-q} + P_p$, $Q = Q_{-p} + Q_q$,
  \marginpar{2=ab}
\begin{equation}
\label{2=ab}
p = mt=ni \geq 2 \text{ and } q=ms=nj \geq 1.
\end{equation}

Suppose that $p=q$, and so $P = P_{-p} + P_p$, $Q = Q_{-p} + Q_p$. Then
$Q=\lambda P_p$ for some $\lambda\in K^*$. Since
$1=[P,Q]=[P, Q-\lambda P]$, $m(P)=2$ and $m(Q-\lambda P)=1$. By the case (iv), the pair
$(P, Q-\lambda P)$ is obtained from the pair $(Y, X)$ by applying an automorphism of the Weyl algebra $A_1$.

So, either $p<q$ or $p>q$. In view of $(P,Q)$-symmetry ($1=[P,Q]=[-Q,P]$), it suffices to
consider, say, the first case only. Since $p<q$, the equalities (\ref{2=ab}) imply that $i<j$ and $t<s$.
Then, using Equation (\ref{degsf}) and the fact that $\deg (p,-p) = p$ for all $p \geq 1$, we see that
$$\deg a = \deg \alpha_{n,i} + \deg \beta'_{m,t} + p -1,$$
$$\deg b = \deg \alpha'_{n,j} + \deg \beta_{m,s} + q -1.$$
Since $i<j$ and $t<s$,  $\deg \alpha_{n,i} < \deg \alpha'_{n,j}$
and $\deg \beta'_{m,t} < \deg \beta_{m,s}$. In particular, $\deg a < \deg b$. This equality
contradicts Equation (\ref{1=ab}) since, by Equation (\ref{degsf}),
$$0=\deg 1 = \deg a - 1 - \deg b +1 = \deg a - \deg b >0.$$
This means that the cases
$p<q$ and $p>q$ are impossible. The proof of the theorem is complete. $\Box$

\begin{corollary}
\label{24Sep16} 
Let $P, Q$ be elements of the first Weyl algebra $A_1$ with $m(P)=1$ or $m(Q)=1$.
If $[P,Q]=1$ then $P = \tau(Y)$ and $Q=\tau(X)$ for some automorphism $\tau\in\Aut_K (A_1)$.
\end{corollary}

{\bf Proof:} Without loss of generality we may assume $m(Q)=1$ and $m(P)\geq 3$. That is $Q = Q_q$ and
$P = \sum_{i\in I} P_i$, where $I \subset \Z$ is a finite set, $q\in \Z\setminus\{0\}$ and the elements $Q_q$ and $P_i$
 are homogeneous in $A_1$. By Equation (\ref{xiaut1}), we may assume that $q>0$.
Then $1 = [P,Q]=\sum_i [P_i, Q_q]$ implies that $-q\in I$, $[P_{-q}, Q_q]=1$ and $[P_j, Q_q]=0$
for all $j\in I$ such that $j\neq -q$.
By Theorem \ref{20Sep16},
$$
q=1, Q_1 = \l X \text{ and } P_{-1} = \l^{-1}Y \text{ for some }\l\in K^*.
$$
By Lemma \ref{a20Sep16}, $C := P - P_{-1} \in C_A(X)=K[X]$. Then $P = \tau(Y)$ and $Q=\tau(X)$
where $\tau: A_1 \to A_1$, $X \mapsto \l X$, $Y \mapsto \l^{-1} Y + C$, is an automorphism.

\section*{Acknowledgements}

This paper has been written during the visit of V.~V.~Bavula to Aachen in 2016,
which was
 supported by the Graduiertenkolleg ``Experimentelle und konstruktive Algebra" of the German Research Foundation (DFG).
The second author has been supported by Project II.6 of SFB-TRR 195 ``Symbolic Tools in Mathematics and their Applications'' of the DFG.

\end{document}